\documentclass[12pt]{article}
\usepackage{amsmath,amsfonts,amssymb,rotating,amsthm, bm}

\usepackage{mathptmx}
\usepackage{times}

\usepackage{enumitem}
\usepackage{caption}
\usepackage{subcaption}
\usepackage{booktabs}
\usepackage {extarrows}
\usepackage{verbatim}
\usepackage{graphicx}
\usepackage{diagbox}
\usepackage{hyperref}

\usepackage{tikz}
\usetikzlibrary{cd}

\allowdisplaybreaks
\def\qed{\nopagebreak\hfill{\rule{4pt}{7pt}}}
\def\proof{\noindent {\it{Proof.} \hskip 2pt}}
\newtheorem{theorem}{Theorem}[section]

\newtheorem{lemma}[theorem]{Lemma} 
\newtheorem{corollary}[theorem]{Corollary}
 
\newtheorem{proposition}[theorem]{Proposition} 

\numberwithin{equation}{section}

\begin{document}

\begin{center}
{\Large\bf Enriched Cycle Structures and 

Roots of Permutations }
 
William Y.C. Chen and Elena L. Wang

Center for Applied Mathematics and KL-AAGDM\\
Tianjin University\\
Tianjin 300072, P.R. China

\vskip 3mm

Emails: { $^1$chenyc@tju.edu.cn, $^{2}$ling\_wang2000@tju.edu.cn}

\vskip 3mm

{\it Dedicated to Persi Diaconis on the occasion of his 80th
birthday}

\end{center}

\begin{abstract} 
This paper is concerned with a duality
  between $r$-regular permutations
and $r$-cycle permutations, and a monotone property 
due to B\'ona-McLennan-White on 
the probability $p_r(n)$ for a random permutation of $\{1,2,\ldots, n\}$
to have an $r$-th root, where $r$ is a prime. For $r=2$, the duality relates permutations with odd cycles to permutations
with even cycles. 
To handle the general case where $r\geq 2$, we define an $r$-enriched permutation as a permutation with 
$r$-singular cycles colored by one of the colors
$1, 2, \ldots, r-1 $.
In this setup, we discover a bijection
between $r$-regular permutations and enriched $r$-cycle permutations, which in turn yields a stronger version of an inequality of B\'ona-McLennan-White. 
This leads to a fully combinatorial understanding of the monotone property, thereby answering their question.   
When $r$ is a prime power $q^l$, we further show that
$p_r(n)$ is monotone.
In the case that
 $n+1 \not\equiv 0 \pmod q$, the equality $p_r(n)=p_r(n+1)$
 has been established by Chernoff.   

\end{abstract}



\usetikzlibrary{decorations.pathreplacing,decorations.markings}

\section{Introduction}

This paper is concerned with
a duality between 
$r$-regular permutations and $r$-cycle permutations,
which are closely related to 
  permutations with an $r$-th root, see, e.g., \cite{Bolker-Gleason-1980, BMW-2000, Kulshammer-etal-2003}.
For an integer $r \ge 2$, a cycle is called $r$-regular if its 
length is not divisible by $r$ and $r$-singular otherwise (i.e., if 
its length is divisible by $r$). Suppose that 
permutations are represented in the cycle notation. A permutation is called $r$-regular provided that all of its cycles are $r$-regular, 
and an $r$-cycle permutation is a permutation where all cycles are $r$-singular. These terms 
were introduced by K\"{u}lshammer, Olsson, and Robinson 
\cite{Kulshammer-etal-2003}. As usual, for $n\geq 1$,
$S_n$ 
stands for 
the set of  permutations of $[n]=\{1,2,\ldots, n\}$. Given a permutation 
$\sigma \in S_n$, it is said to have 
an $r$-th root if there exists a permutation $\pi \in S_n$ 
such that  $\pi^r = \sigma$. 
Permutations with an $r$-th root can be characterized
in terms of their cycle types \cite[p. 158]{Wilf-2005}.   
For the number of such permutations, one may refer to the sequence A247005 in the OEIS \cite{OEIS}. The exponential generating function for this count was discussed in  Bender \cite{Bender-1974}, 
see also \cite[p. 159]{Wilf-2005}. 

We shall follow the terminology in 
\cite{Kulshammer-etal-2003}. Throughout the paper, 
 ${\rm Reg}_r(n)$ and ${\rm Cyc}_r(n)$ will
 stand for the set of $r$-regular permutations of $[n]$ and
 the set of $r$-cycle permutations of $[n]$, respectively. 
 Note that 
 ${\rm NODIV}_r(n)$ and ${\rm PERM}_r(n)$ are used in 
 \cite{BMW-2000}.
  For $r\geq 2$ and $n=0$, set $\left\lvert{\rm Reg}_r(0)\right\rvert=\left\lvert{\rm Cyc}_r(0)\right\rvert=1$.
 Clearly, $\left\lvert{\rm Cyc}_r(n)\right\rvert=0$ whenever $n \not\equiv 0 \pmod r$.

 The enumeration of   $r$-regular permutations  
dates back to 
  Erd\H{o}s and Tur\'an \cite{Erdos-Turan-1967}.
  By using generating functions, they
  showed that for $n\geq 1$, and $r$ a prime power,
  the proportion of $r$-regular 
  permutations in $S_n$ 
  equals 
\[ \prod_{k=1}^{\lfloor n/r \rfloor}\frac{rk-1}{rk}.\]
It was later observed that 
the above formula remains true for an arbitrary integer  $r \ge 2$, for example, see \cite{Maroti-2007}.

There are various ways to count ${\rm Reg}_r(n)$ and ${\rm Cyc}_r(n)$, 
see \cite{Beals-etal-2002, Bertram-Gordon-1989, 
Bolker-Gleason-1980, BMW-2000, Glasby-2001, Maroti-2007, Wilf-2005}. 
In particular, for $r\geq 2$,
B\'ona, McLennan and White \cite{BMW-2000}
presented a bijective argument 
to deduce the number of $r$-regular permutations of $[n]$
from the number of $r$-regular permutations of $[n-1]$. 
As a consequence, they confirmed the conjecture of  Wilf \cite{Wilf-2005} that 
the probability $p_2(n)$ for a random permutation of $[n]$ to have
a square root is monotonically nonincreasing in $n$. 
Such permutations have been called square permutations \cite{Blum-1974}. For example, 
 $(1\ 2\ 3\ 4)\,(5\ 6\ 7\ 8)$ is a square permutation and it has a square root $(1\ 5\ 2\ 6\ 3\ 7\ 4\ 8)$. 
Beyond the case $r=2$, they proved that, more generally, for any prime $r$, the probability 
$p_r(n)$ that a random permutation of $[n]$ has an $r$-th root is nonincreasing in $n$.

Notice that the monotone property does not hold in general.
 For example, when $r=6$, we have $p_6(4)=1/6$ but $p_6(5)=1/3$.
Nonetheless, B\'ona, McLennan and White \cite{BMW-2000} showed 
that for any $r\geq 2$,
\begin{align*}
    p_r(n) \to 0,
\end{align*}
as $n \to \infty$.

Table \ref{tb-1} exhibits the values  of $p_r(n)$ for $r=2,3,5$ and 
$1\leq n \leq 12$. 
\begin{table}[ht]
        \setlength{\tabcolsep}{8pt}
        \renewcommand{\arraystretch}{1.2}
	\centering
	\begin{tabular}{|c|cccccccccccc|}
		\hline
           \diagbox[width=2.8em,height=2.5em]{$r$}{\raisebox{4pt}{$n$}}
            & \raisebox{-1pt}{$1$} & \raisebox{-1pt}{$2$} & \raisebox{-1pt}{$3$} & \raisebox{-1pt}{$4$} & \raisebox{-1pt}{$5$} & \raisebox{-1pt}{$6$} & \raisebox{-1pt}{$7$} & \raisebox{-1pt}{$8$} & \raisebox{-1pt}{$9$} &\raisebox{-1pt}{$10$} & \raisebox{-1pt}{$11$} & \raisebox{-1pt}{$12$} \\	\hline
		$2$ & ${ 1}$ & ${ \dfrac{1}{2}}$ & ${ \dfrac{1}{2}}$ & ${ \dfrac{1}{2}}$ & ${ \dfrac{1}{2}}$ & ${ \dfrac{3}{8}}$ & ${ \dfrac{3}{8}}$ & ${ \dfrac{17}{48}}$ & ${ \dfrac{17}{48}}$ & ${ \dfrac{29}{96}}$ & ${ \dfrac{29}{96}}$ & ${ \dfrac{209}{720}}$ \raisebox{23pt}{}\\[10pt]	\hline
		$3$ &  ${ 1}$ & ${ 1}$ & ${ \dfrac{2}{3}}$ & ${ \dfrac{2}{3}}$ & ${ \dfrac{2}{3}}$ & ${ \dfrac{5}{9}}$ & ${ \dfrac{5}{9}}$ & ${ \dfrac{5}{9}}$ & ${ \dfrac{1}{2}}$ & ${ \dfrac{1}{2}}$ & ${ \dfrac{1}{2}}$ & ${ \dfrac{37}{81}}$ \raisebox{23pt}{}\\[10pt]		\hline
        $5$ & ${ 1}$ & ${ 1}$ & ${ 1}$ & ${ 1}$ & ${ \dfrac{4}{5}}$ & ${ \dfrac{4}{5}}$ & ${ \dfrac{4}{5}}$ & ${ \dfrac{4}{5}}$ & ${ \dfrac{4}{5}}$ & ${ \dfrac{18}{25}}$ & ${ \dfrac{18}{25}}$ & ${ \dfrac{18} {25}}$  \raisebox{23pt}{} \\[10pt]		\hline
	\end{tabular}
        \caption{The values of $p_r(n)$.}
        \label{tb-1}
\end{table}

As set forth by B\'ona, McLennan and White, 
their proof of the 
monotone property is mostly  combinatorial, and they 
left the question of seeking a 
fully combinatorial reasoning, 
which amounts to a combinatorial understanding of 
the following inequality 
\begin{align} \label{BMW-I}
    \left\lvert\mbox{Cyc}_{r^2}(mr^2) \right\rvert \le \left\lvert \mbox{Reg}_r(mr^2)\right\rvert, 
\end{align}
which we shall  call  
the B\'ona-McLennan-White inequality, 
or the BMW inequality, for short.

Note that for $r=2$,  ${\rm Reg}_2(n)$ and ${\rm Cyc}_2(n)$  are usually written as
${\rm Odd}(n)$ and ${\rm Even}(n)$, respectively. 
In the literature, $2$-regular permutations
are also known as odd order permutations, which
are related to ballot permutations, see, for example,
\cite{Bernardi-etal-2010, Lin-Wang-Zhao-2022, Spiro-2020}. However, even order permutations
are referred as permutations with at least one even
cycle. These terms originated from the 
notion of the order of an element in a group.

It is a known result that ${\rm Odd}(2n)$ and ${\rm Even}(2n)$ share the same cardinality, $((2n-1)!!)^2$, 
see A001818 in the OEIS \cite{OEIS}. The equality of their cardinalities entails the existence of a bijection between ${\rm Odd}(2n)$ and ${\rm Even}(2n)$; however, this mapping is not immediately evident upon first inspection.
A specific correspondence was found by Sayag based on the canonical representation of permutations, see
  B\'ona \cite[Lemma 6.20]{Bona-2005}. 
  An intermediate structure, 
  which we call nearly odd order permutations,
  was introduced in \cite{Chen-2024}. It induces incremental transforms from a
  permutation in ${\rm Odd}(2n)$ to a permutation in 
  ${\rm Even}(2n)$.  

It is natural to ask whether the correspondence between ${\rm Odd}(2n)$ and ${\rm Even}(2n)$ can be extended to any  $r \ge 2$? 
In this paper, we introduce the structure of enriched permutations, 
and we find a bijection
between  $r$-regular permutations of $[rn]$ 
and enriched $r$-cycle permutations of $[rn]$. 
As an immediate consequence, we achieve a combinatorial
 comprehension
 of the BMW inequality, 
 or a stronger version, strictly speaking. 
This answers the question of
B\'ona, McLennan and White for any prime $r \geq 3$.
As for the case $r=2$, we fill up with some discussions 
for the sake of completeness. 

While the monotone property
of $p_r(n)$ does not hold for
general $r$, we show that 
it is valid for prime powers $r=q^l$. The proof relies on
 a stronger version of the 
B\'ona-McLennan-White inequality and the characterization 
of permutations with an $r$-th root, for a prime power $r$,
due to A. Knopfmacher and R. Warlimont \cite[p. 158]{Wilf-2005}.

\section{$r$-Enriched permutations}

\label{Bijective proof}

The aim of this section is to provide
a bijection between $r$-regular permutations of $[rn]$  and
enriched $r$-singular permutations of $[rn]$, for $r\geq 2$ and
$n\geq 1$. Given $r \ge 2$, by saying that a permutation is $r$-enriched we mean
that each $r$-singular cycle is colored by one of the
$r-1$ colors in the set \(\{1, 2, \dots, r-1\}\). Later, we simply call such a permutation enriched.
Bear in mind that $r$-regular cycles are never colored. 

Given $r \ge 2 $, we shall use the symbol $*$ to 
 signify the enriched structure. For example, ${\rm Cyc}^*_r(rn)$ 
denotes the set of enriched $r$-cycle permutations 
of $[rn]$. 
Throughout, we represent each permutation in cycle notation, where each cycle is written as a linear order starting with its smallest element, and the cycles are arranged in increasing order of their minima. 
We use the subscript of a cycle to denote the color assigned to it. 
For example, for $r=3$,    
$(1\ 2\ 4)_2\,(3)\,(5\ 6)$ represents an $3$-enriched
permutation for which  the $3$-singular cycle $(1\ 2\ 4)$ 
is colored by $2$.

To transform an $r$-regular permutation of $[rn]$ to an enriched $r$-singular permutation of $[rn]$, we introduce an intermediate structure like nearly odd permutations emerging
in \cite{Chen-2024}. For $n\geq 1$, we say that a permutation $\sigma$ of $[n]$
is nearly $r$-regular if its cycles are all $r$-regular 
except that the one containing $1$ is $r$-singular. 
The notation ${\rm NReg}_r(n)$ stands for
the set of all nearly $r$-regular permutations of $[n]$. 
For example, $(1\ 2\ 4)\,(3)\,(5\ 6)$ is a nearly $3$-regular permutation. 

Enriched nearly $r$-regular permutations serve as an intermediate structure linking $r$-regular and enriched $r$-cycle permutations. More
precisely, we manipulate the cycle structures to construct a bijection between ${\rm Reg}_r(rn)$ 
and ${\rm Cyc}^*_r(rn)$. 
For $r=2$, it reduces to a bijection between ${\rm Odd}(2n)$ and ${\rm Even}(2n)$.

\begin{theorem}
    \label{Phi}
    For any $r \ge 2$, there is a bijection $\Phi$ from ${\rm Reg}_r(rn)$ 
    to ${\rm Cyc}^*_r(rn)$. Moreover, if $\sigma \in {\rm Reg}_r(rn)$ 
    and the cycle containing $1$ in $\sigma$ has 
    length $l=rk+i$, $1 \le i \le r-1$, then $\Phi(\sigma) \in {\rm Cyc}^*_r(rn)$, 
    where the cycle containing $1$ in $\Phi(\sigma)$ has length $rk+r$.
\end{theorem}

To prove the theorem, let $Q_{r,\,k}(n)$ denote the set of permutations of $[n]$ for which  
  the length of the cycle containing $1$ is $k$, 
and the other cycles are $r$-regular. 
We first build a bijection between 
$Q_{r,\,k}(n)$ and $Q_{r,\,k+1}(n)$ by applying  
an elegant bijection of B\'ona, McLennan and White 
in \cite[Lemma 2.1]{BMW-2000},
which is a paradigm of a 
recursive algorithm. 

\begin{lemma}
    \label{BMW-bij-1}
    For all $r \geq 2$ and $n+1 \not\equiv 0 \pmod r$, there is 
    a bijection $\Psi$ from ${\rm Reg}_r(n) \times[n+1]$ to ${\rm Reg}_r(n+1)$. 
\end{lemma}

Practically, to make use of the bijection, it is unnecessary to adjust the elements of the underlying sets to
fit in the above canonical form. It seems to be  
convenient to harness the following more generic 
formulation, and it might be informative to
reproduce the proof as such. For a nonempty set $S$, we use 
${\rm Reg}_r(S)$ to denote the
set of $r$-regular permutations of $S$. 

\begin{lemma}
\label{BMW-bij-2}
Let $S$ be a nonempty finite set and $r$ an integer such that $r\geq 2$. If $\lvert S \rvert\not \equiv 0 \pmod{r}$, then
there is a bijection $\Delta$ from ${\rm Reg}_r(S)$ to the 
set of pairs $(x, \pi)$ such that 
 $x\in S$ and $\pi$ is in ${\rm Reg}_r(S
 \setminus \{x\})$.   
\end{lemma}

\proof
Assume that $\sigma$ belongs to ${\rm Reg}_r(S)$ and  
$\lvert S \rvert  \not\equiv 0 \pmod r$. 
From now on, we shall use $\lvert \sigma \rvert$ 
to denote the number of elements of $\sigma$. 
Let $D_1$ denote the first cycle of $\sigma$, $l$ be its length, $x$ be the last entry in $D_1$, and let 
$\Tilde{\sigma}$ be $\sigma$ with $D_1$ being removed. 
In effect,
the map $\Delta$ will remove
$x$ in $D_1$  and turn it into
a distinguished element. 
We encounter three cases. 

\noindent
Case 1: $l=1$. Then set $\Delta(\sigma)=(x, \Tilde{\sigma})$. By construction, the element $x$ is smaller than every element of $\Tilde{\sigma}$.

\noindent
Case 2: $l \not\equiv 1 \pmod r$. Then
remove $x$ from $D_1$ to get $C_1$ and set $\Delta(\sigma)=(x, C_1\,\Tilde{\sigma})$. 
Contrary to the previous case, the element $x$ is greater than the smallest 
element of $C_1\,\Tilde{\sigma}$. 
Since $l \not\equiv 0, 1 \pmod r$, we have
$\left\lvert C_1 \right\rvert = l-1 \not\equiv -1,0 \pmod r$, which ensures that the resulting permutation is $r$-regular.

\noindent
Case 3: $l \equiv 1 \pmod r$ and $l \neq 1$. 
Let $\Tilde{x}$ be the second-to-last element in $D_1$ 
and $C_1$ be the cycle obtained
from $D_1$ by removing $x$ and $\Tilde{x}$. Since
$\left\lvert \Tilde{\sigma}\right\rvert+1
=\lvert \sigma \rvert-l+1 \not\equiv 0 \pmod r$, we can 
apply $\Delta^{-1}$ to $(\Tilde{x}, \Tilde{\sigma})$ to get $\Tilde{\pi}$. Then set
 $\Delta(\sigma)=(x, C_1\, \Tilde{\pi})$. 
In this situation, the element $x$ is greater than the smallest 
element of $C_1\, \Tilde{\pi}$ and $\left\lvert C_1\right\rvert = l-2 \equiv -1 \pmod r$.

It remains to verify that $\Delta$ is a bijection. 
Given a pair $(x,\pi)$ where $x\in S$ and 
$\pi$ is an $r$-regular permutation of $S 
\setminus \{x\}$ with $\lvert \pi \rvert+1
\not\equiv 0 \pmod r$. 
Let $C_1$ denote the first cycle of $\pi$, 
$l$ be its length and let $\Tilde{\pi}$ be permutation $\pi$ with $C_1$ being removed. 
Conversely, the map $\Delta^{-1}$ will place $x$ as the last entry 
in the first cycle of $\Delta^{-1}(x,\pi)$.  
Accordingly, we face three possibilities.

\noindent
Case 1: The element $x$ is smaller than every element of $\pi$.  
Then set $\Delta^{-1}(x,\pi)=(x) \, \pi$. In this case, 
$\left\lvert D_1\right\rvert$=1.

\noindent
Case 2: The element $x$ is greater than the smallest 
element of $\pi$ and $l \not\equiv -1 \pmod r$.
Let $D_1$ be $C_1$ with $x$ appended 
to the end of $C_1$. Then set 
$\Delta^{-1}(x,\pi)=D_1 \, \Tilde{\pi}$. 
Notice that $l \not\equiv -1, 0 \pmod r$, 
and so $\left\lvert D_1 \right\rvert=l+1 \not\equiv 0, 1 \pmod r$. 

\noindent
Case 3: The element $x$ is greater than the smallest 
element of $\pi$ and $l \equiv -1 \pmod r$.
Under the conditions 
$\lvert \pi \rvert-l \equiv \lvert \pi \rvert+1 \pmod r$ and $\lvert \pi \rvert +1 \not\equiv 0 \pmod r$, we  have $\left\lvert \Tilde{\pi} \right\rvert=\lvert \pi \rvert-l \not\equiv 0 \pmod r$. 
Thus we can apply $\Delta$ to $\Tilde{\pi}$ to get $(\Tilde{x}, \Tilde{\sigma})$. Let $D_1$ 
be $C_1$ with $\Tilde{x}x$ appended
to the end.
Then set $\Delta^{-1}(x,\pi)=D_1 \, \Tilde{\sigma}$. 
Notice that in this case $\left\lvert D_1 \right\rvert = l +2 \equiv 1 \pmod r$ and $\left\lvert D_1 \right\rvert \neq 1$, completing the proof. 
\qed

For example, when $r=3$, 
consider \[ \sigma = (1\ 8\ 2\ 5)\,(3)\,(4)\,(6\ 7).\] Since the length of the first cycle is congruent to $1 \pmod 3$, we are in Case $3$. 
Hence,
\begin{align*}
    \Delta(\sigma)=(5,(1\ 8)\,\Tilde{\pi}),
\end{align*}
where 
\begin{align*}
    \Tilde{\pi}=\Delta^{-1}(2,(3)\,(4)\,(6\ 7)).
\end{align*}
Because the element $2$ is smaller than every element in $\Tilde{\pi}$, we are in Case $1$ of $\Delta^{-1}$. Therefore, 
\begin{align*}
    \Delta^{-1}\left(2,(3)\,(4)\,(6\ 7)\right)=(2)\,(3)\,(4)\,(6\ 7).
\end{align*}
Consequently,
\begin{align*}
     \Delta(\sigma)=(5,(1\ 8)\,(2)\,(3)\,(4)\,(6\ 7)).
\end{align*}

We now turn to the description of the bijection $\varphi$.  
The following lemma is the building block of the
correspondence between $r$-regular permutations and 
enriched $r$-cycle permutations. It rests on
the Lemma of B\'ona, McLennan and White, 
as restated in Lemma \ref{BMW-bij-2}. 

\begin{lemma}
\label{varphi}
    Let $n, r, k$ be integers such that
    $n \ge 0, r \ge 2, k \ge 1$ and $n -k \not\equiv 0 \pmod r$. Then there is a bijection $\varphi$ 
    from $Q_{r,\,k}(n)$ to $Q_{r,\,k+1}(n)$.
\end{lemma}

\proof
We proceed to construct a
map $\varphi$ from $Q_{r,\,k}(n)$ to $Q_{r,\,k+1}(n)$ by
employing the above bijection $\Delta$.  
Assume that
$n -k \not\equiv 0 \pmod r$ and $k \ge 1$. 
Let $\sigma \in Q_{r,\,k}(n)$, and denote by $C_1$ the cycle containing $1$. Define  $\Tilde{\sigma} = \sigma - C_1$, by which we mean the permutation obtained from $\sigma$ by removing $C_1$.
Since $\left\lvert\Tilde{\sigma}\right\rvert=n-k \not\equiv 0 \pmod r$ and 
$\Tilde{\sigma}$ is an $r$-regular permutation, 
applying the map of $\Delta$,
we get $\Delta(\Tilde{\sigma})=(x, \Tilde{\pi})$.
Now, let $D$ denote the cycle obtained 
  from $C_1$ with $x$ attached at the end. 
  Set 
  $\varphi (\sigma) = D\, \Tilde{\pi}$, which is readily seen to lie in $Q_{r,\,k+1}(n)$. 

Conversely, let us define a map $\alpha$ from $Q_{r,\,k+1}(n)$ to
$Q_{r,\,k}(n)$. 
Given $\pi \in Q_{r,\,k+1}(n)$, 
where $n-k \not\equiv 0 \pmod r$ and $k \ge 1$, 
let $C_1$ be the first cycle of $\pi$, and 
let $D$ be the cycle obtained from $C_1$
by removing its last entry $x$. 
Define $\Tilde{\pi}$ be $\pi$ with $C_1$ being removed.
Note that $\left\lvert \Tilde{\pi} \right\rvert +1=n-k  \not\equiv 0 \pmod r$ and
$ \Tilde{\pi}$ is an $r$-regular permutation. 
Then set
\[ \alpha(\pi) = D \, \Delta^{-1}(x, \Tilde{\pi}),\]
which belongs to $Q_{r,\,k}(n)$.

It is straightforward to verify that the maps
$\varphi$ and $\alpha$ are well-defined and are inverses of each other. Thus
$\varphi$ is a bijection. 
\qed

For example, if $r=3$, 
then $\varphi((3)\,(5\ 6))=(3\ 6)\,(5)$ and $\varphi((3 \ 6)\,(5))=(3\ 6\ 5)$.

   Writing $n-k=mr+d$ with $0 < d < r $, 
 it is known that, see \cite{Beals-etal-2002, Bertram-Gordon-1989, 
Bolker-Gleason-1980, BMW-2000, Glasby-2001, Maroti-2007, Wilf-2005}, 
 \begin{align*}
   \left\lvert {\rm Reg}_{r}(n-k)\right\rvert =   (n-k)! \, \frac{(r-1)\, (2r-1)\cdots(mr-1)}{r^m m!},
\end{align*}
from which we deduce that 
 \begin{align}
   \left\lvert Q_{r,\,k}(n)\right\rvert =   \left\lvert Q_{r,\,k+1}(n)\right\rvert =   (n-1)! \, \frac{(r-1)\, (2r-1)\cdots(mr-1)}{r^m m!}.
\end{align} 

For $k \ge 1$, let $A_{n,\,2k-1}$ denote 
the set of permutations of $[n]$ with only odd cycles for which the
element $1$ appears in a cycle of length $2k-1$, and
let $P_{n,\,2k}$ denote the set of permutations of $[n]$ with 
odd cycles except that the element $1$ is contained in an even 
cycle of length $2k$. When $r=2$, we come to the following correspondence.
Notice that the 
 construction in \cite{Chen-2024} by way of 
 breaking cycles does not possess this refined property.

\begin{corollary}
    \label{bij-r-2}
     For $k \ge 1$, the map $\varphi$ defined in Lemma \ref{varphi} gives a bijection between  $A_{2n,\, 2k-1}$ and $P_{2n, \,2k}$, as well as a bijection between  $P_{2n+1,\, 2k}$ and $A_{2n+1,\, 2k+1}$.
\end{corollary}

For example, when $r=2$, 
given $\sigma = (1\ 2\ 3\ 4\ 6)\,(5\ 10\ 8)\,(7)\,(9) \in A_{10,\,5}$, we have
  \begin{align*}
      \Delta((5\ 10\ 8)\,(7)\,(9))=(8,(5)\,(7\ 9\ 10)).
  \end{align*}
Thus
\begin{align*}
    \varphi(\sigma) & = (1\ 2\ 3\ 4\ 6\ 8) \,(5) \, (7\ 9\ 10) \in P_{10,\,6}.
\end{align*}

Below are the explicit formulas:
\begin{align}
    & \left\lvert A_{2n,\,2k-1} \right\rvert =  \left\lvert P_{2n,\,2k}\right\rvert=\frac{(2n-1)!}{(2n-2k)!}\, \left((2n-2k-1))!!\right)^2,
    \\[6pt]
    & \left\lvert P_{2n+1,\,2k}\right\rvert =  \left\lvert A_{2n+1,\,2k+1} \right\rvert=\frac{(2n)!}{(2n-2k)!}\, \left((2n-2k-1))!!\right)^2.
\end{align}

Exploiting the bijection $\varphi$, we are led to the following incremental 
transformation $\Lambda$. 

\begin{theorem}
    \label{Lambda}
        For all $r \ge 2$, there is a bijection $\Lambda$ 
        from ${\rm Reg}_r(rn)$ to ${\rm NReg}^*_r(rn)$. 
        Moreover, if $\sigma \in {\rm Reg}_r(rn)$ and the 
        cycle containing $1$ in $\sigma$ has length $l=rk+i$, 
        $1 \le i \le r-1$, then $\Lambda(\sigma) \in {\rm NReg}^*_r(rn)$, 
        where the cycle containing $1$ in $\Lambda(\sigma)$ has length $rk+r$.
\end{theorem}

\proof
Let $\sigma$ in ${\rm Reg}_r(rn)$. Assume 
that its first cycle length is $rk+i$, where $1 \le i \le r-1$. 
Since $rn \not\equiv rk+i \pmod r$, 
we can apply the bijection $\varphi$ in Lemma
\ref{varphi} to $\sigma$ 
to get a permutation $\pi$. There are two
possibilities with regard to the length of the
first cycle of $\pi$.

If $\pi$ is in ${\rm NReg}_r(rn)$, in which case, $i+1=r$, we color its first cycle with $r-1$ to obtain $\Lambda(\sigma) \in {\rm NReg}^*_r(rn)$.  

If $\pi$ stays in ${\rm Reg}_r(rn)$ with
the length of the first cycle increased by $1$ 
in comparison with $\sigma$,  that is, 
the length of the first cycle of $\pi$ equals $rk+i+1$
with $rk+i+1 \not\equiv 0 \pmod r$. 
Again, since $rn \not\equiv rk+i+1 \pmod r$,
we may move on to apply the bijection $\varphi$ once more. 
The procedure continues until we obtain a permutation 
$\pi$ in ${\rm NReg}_r(rn)$. Color $\pi$'s first cycle with $i$ and define the resulting enriched permutation to be $\Lambda(\sigma)$, which lies in ${\rm NReg}^*_r(rn)$. Apparently, it takes
$r-i$ steps to reach this point.

As an example, for $r=3$, 
we have \[\Lambda\left((3)\,(5\ 6)\right)=(3\ 6\ 5)_1\]
and 
\[\Lambda\left((1\ 2)\,(3\ 4)\,(5\ 6)\right)  =  (1\ 2\ 4)_2\,(3)\,(5\ 6). \]
 It is readily 
seen that the process is reversible because
the color of the $r$-singular cycle keeps track of
the number of times that the map $\varphi$ is applied. 
Thus $\Lambda$ is a bijection.
\qed

\noindent \textit{Proof of Theorem \ref{Phi}.} 
We proceed to define a map $\Phi$ from  ${\rm Reg}_r(rn)$ to  ${\rm Cyc}_r^*(rn)$  by successively applying the map $\Lambda$ in Theorem \ref{Lambda}.

Given an $r$-regular permutation $\sigma$, our goal is to create a sequence of colored $r$-singular cycles $C_i^*$. At the first step, applying the bijection $\Lambda$ to $\sigma$, we get 
\[\Lambda(\sigma)=C_1^*\, \sigma_1,\] where $C_1^*$ is a colored $r$-singular cycle and $\sigma_1$ is an $r$-regular permutation (possibly empty), and $C_1^* \sigma_1$ stands for the
permutation obtained by
putting together the colored cycle $C_1^*$ and the cycles in $\sigma_1$.   If $\sigma_1$ is empty, then we set $\Phi(\sigma)=C_1^*$; otherwise,  applying $\Lambda$ again to $\sigma_1$  gives \[\Lambda(\sigma_1)=C_2^*\, \sigma_2.\] If $\sigma_2$ is empty, we define $\Phi(\sigma)=C_1^*\, C_2^*$; otherwise, we may apply $\Lambda$ to $\sigma_2$.  

We may iterate the procedure
as follows. Assume that we have
obtained a nonempty permutation $\sigma_i$ for some $i$.  Applying $\Lambda$ to $\sigma_i$ yields 
\[\Lambda(\sigma_i)=C_{i+1}^* \sigma_{i+1},\] 
where $C_{i+1}^*$ is a colored $r$-singular cycle and $\sigma_{i+1}$ is an $r$-regular permutation (possibly empty).
If $\sigma_{i+1}$ is not empty, we apply $\Lambda$ again; otherwise, we stop.  Clearly, the above process will terminate at some point when $\sigma_{i+1}$ is empty. Now define 
\[\Phi(\sigma)=C_1^*\, \cdots \, C_{i+1}^*.\]
Since $\Lambda$ is bijective, so is $\Phi$.\qed 

For example, for $r=3$, given 
\[\sigma=(1\ 2)\,(3\ 4)\,(5\ 6).\]
Applying  $\Lambda$ to $\sigma$ yields 
\[\Lambda \left((1\ 2)\,(3\ 4)\,(5\ 6)\right) =  (1\ 2\ 4)_2\,(3)\,(5\ 6);\] at this stage, the permutation still contains $r$-regular cycles \((3)\) and \((5\ 6)\), 
so we need to repeat the above procedure. One more round of 
iteration gives
\[\Lambda \left((3)\,(5\ 6)\right) =  (3\ 6\ 5)_1,\] 
thus we finally obtain
 \[
    \Phi\left((1\ 2)\,(3\ 4)\,(5\ 6)\right)=  (1\ 2\ 4)_2\,(3\ 6\ 5)_1.
\]

\section{The B\'ona-McLennan-White inequality }
\label{BMW-application}

In the proof of the following monotone property, 
there is an inequality that demands a
combinatorial explanation. As will be seen,
the structure of enriched cycle permutations entails a combinatorial interpretation of this inequality. Recall that $p_r(n)$ is
the probability for a random permutation of $[n]$ to have
an $r$-th root.

\begin{theorem}[B\'ona, McLennan and White \cite{BMW-2000}]
\label{main-BMW}
    For all positive integers $n$ and all primes $r$,  we have, 
    \begin{align*}
        p_r(n) \ge p_r(n+1).
    \end{align*}
\end{theorem}

The above assertion consists  of three 
circumstances contingent to modulo conditions on $n+1$.

\begin{theorem}[B\'ona, McLennan and White \cite{BMW-2000}]
\label{root-probability}
Let $r$ be a prime. Then we have the following. 
\begin{enumerate}
    \item[(i) ] If $n+1 \not\equiv 0 \pmod r$, then $p_r(n)=p_r(n+1)$.
    \item[(ii) ] If $n+1 \equiv 0 \pmod {r}$ but $n+1 \not\equiv 0 \pmod {r^{2}}$, 
    then 
    \begin{align*}
        p_r(n) \geq \frac{n+1}{n}p_r(n+1)
    \end{align*} 
    with equality only 
    when $n+1=kr$, where $k=1,2,\ldots,r-1$.
    \item[(iii) ] If $n+1 \equiv 0 \pmod {r^{2}}$, then $p_r(n) \ge p_r(n+1)$ with 
    equality only when $r=2$ and $n=3$.
\end{enumerate}
\end{theorem}

The proof of the above theorem builds  
upon a special case of the characterization of permutations with an 
$r$-th root, due to Knopfmacher and Warlimont, see
  \cite[p. 158]{Wilf-2005}. In particular, when $r$ is prime,  a 
permutation has an $r$-th root if and only if  for 
any positive integer $i$,
the number of cycles of length $ir$ is a multiple of $r$. 
Making use of the bijection $\Psi$ as restated in Lemma \ref{BMW-bij-2}, 
B\'ona, McLennan and White gave an entirely combinatorial proof of (i) and (ii).  However, in order to have a fully combinatorial 
understanding of (iii), one
needs a combinatorial interpretation of the following inequality, 
see \cite[Lemma 3.3]{BMW-2000}, which we call 
the B\'ona-McLennan-White inequality, or the BMW inequality, for short. 

\begin{lemma}
    For all $r \geq 2$ and $m \ge 1$, 
     \begin{align}
     \label{BMW-inequality}
    \left\lvert {\rm Cyc}_{r^2}(mr^2) \right\rvert < 
    \left\lvert {\rm Reg}_r(mr^2)\right\rvert .
    \end{align}
\end{lemma}
 
  The BMW-inequality was proved in \cite{BMW-2000} 
 by means of generating functions. 
 In fact,  we observe that a stronger version of 
 \eqref{BMW-inequality}, i.e., Theorem \ref{BMW-new}, holds, which can be deduced 
 from the following known formulas, 
 see \cite{Beals-etal-2002, Bertram-Gordon-1989, 
Bolker-Gleason-1980, BMW-2000, Glasby-2001, Maroti-2007, Wilf-2005}.
 For $r \ge 2$ and $m \ge 1$, 
\begin{align}
    \label{Formula-Cyc}
    \left\lvert{\rm Cyc}_{r}(rm)\right\rvert& =(rm)!\frac{ (1+r)\,(1+2r)\cdots(1+(m-1)r)}{r^m m!},
    \\[6pt]
    \label{Formula-Reg}
    \left\lvert{\rm Reg}_{r}(rm)\right\rvert& =(rm)!\frac{(r-1)\, (2r-1)\cdots(mr-1)}{r^m m!}.
\end{align} 
On the other hand, it is transparent from a combinatorial
point of view.

\begin{theorem}
    For  $r \geq 2$ and $n\geq 1$, 
     \begin{align}
     \label{BMW-new}
    \left\lvert{\rm Cyc}_r(n)\right\rvert \le \left\lvert{\rm Reg}_r(n)\right\rvert,
    \end{align}
   where the equality holds only when
     $r=2$ and $n$ is even. 
\end{theorem}

\proof
When $n \not\equiv 0 \pmod r$, we have $\left\lvert{\rm Cyc}_{r}(n)\right\rvert=0$, 
  nothing needs to be done. When $n = rm$, 
by restricting to only one color, we see
that
\begin{equation}
      \left\lvert{\rm Cyc}_r(rm)\right\rvert \leq  \left\lvert{\rm Cyc}^*_r(rm)\right\rvert.
\end{equation}
But Theorem \ref{Phi} says that $ \left\lvert{\rm Reg}_{r}(rm)\right\rvert=  \left\lvert{\rm Cyc}^*_r(rm)\right\rvert$, and so  \eqref{BMW-new} follows. The equality holds only when $r=2$ and $n$ is even. This concludes the proof.
\qed

To see that the BMW inequality \eqref{BMW-inequality} 
stems from \eqref{BMW-new}, just observe that for $r\geq 2$,
 \begin{align*}
  {\rm Cyc}_{r^2}(mr^2) \subset {\rm Cyc}_r(mr^2).
\end{align*}
This inequality, together with 
   the combinatorial reasoning in \cite{BMW-2000} 
   gives rise to
   the 
   conclusion that $p_r(n) > p_r(n+1)$ for any prime  $r \ge 3$ and $n+1 \equiv  0 \pmod {r^2}$.

As defined before,
$\left\lvert{\rm Reg}_r(0)\right\rvert=1$ and $\left\lvert{\rm Cyc}_r(0)\right\rvert=1$. 
With Lemma \ref{BMW-bij-1} in hand, 
it is easy to get the following recurrence of 
$\left\lvert{\rm Reg}_r(rm)\right\rvert$. For details, we refer to Lemma 2.1 and Lemma 2.6 in \cite{BMW-2000}. 

\begin{lemma}
\label{recurrence-reg}
    For all $r \ge 2$ and $m \ge 1$, we have 
    \begin{align}
        &\left\lvert{\rm Reg}_r(rm)\right\rvert = (rm-1)\, (rm-1)_{r-1} \left\lvert{\rm Reg}_r(rm-r)\right\rvert, 
         \label{r-1}
    \end{align} 
    where $(x)_{m}$ stands for the lower factorial $x(x-1)\cdots(x-m+1)$. 
\end{lemma}

The above  relation can also be
deduced inductively by using the recursive 
generation of permutations in the cycle notation, 
see, for example, \cite{Beals-etal-2002, Herrera-1957}. 
In \cite{Beals-etal-2002}, it has been shown that‌ 
\begin{align*}
         \begin{split}\left\lvert {\rm Reg}_r(rm)\right\rvert &=   \sum_{1 \le l \le r-1} (rm-1)_{l-1}\left\lvert {\rm Reg}_r(rm-l)\right\rvert \\[3pt]
         & \qquad + (rm-1)_{r}  \left\lvert {\rm Reg}_r(rm-r)\right\rvert, 
         \end{split}
\end{align*} 
and by an easy induction on $n$, it can be shown that  
\begin{align*}
    \left\lvert{\rm Reg}_r(rm-l)\right\rvert=(rm-l)_{r-l} \left\lvert{\rm Reg}_r(rm-r)\right\rvert, \,1 \le l \le r-1.
\end{align*}
This proves \eqref{r-1}.

Similarly, we have the following recurrence relation for
${\rm Cyc}_r(rm)$.

\begin{lemma}
\label{recurrence-cyc}
    For all $r \ge 2$ and $m \ge 1$, we have 
    \begin{align}
        &\left\lvert{\rm Cyc}_r(rm)\right\rvert = (rm-1)_{r-1}(rm-r+1) \left\lvert{\rm Cyc}_r(rm-r)\right\rvert. \label{r-2}
    \end{align}
\end{lemma}

\proof
Let $\sigma$ be a permutation  in ${\rm Cyc}_r(rm)$.
Let $l$ be the length of the first cycle of 
$\sigma$. If $l = r$, then
there are $(rm-1)_{r-1}$ choices to form the first cycle. 
If the first cycle contains 
more than $r$ elements, say $(1\ \cdots\ j_r\ j_{r+1}\ \cdots)$, then there are $(rm-1)_{r}\left\lvert{\rm Cyc}_r(rm-r)\right\rvert$ choices. We can break the first
cycle into two segments $ 1\ \cdots\ j_r $
and $j_{r+1}\ \cdots$. The second segment can
be viewed as a cycle with a distinguished element
$j_{r+1}$. Combining this cycle with a distinguished
element and other cycles, we see a 
permutation in ${\rm Cyc}_r(rm-r)$ with a distinguished
element. There are  $(rm-1)_{r-1}$ for the
first segment $ 1\ \cdots\ j_r $ and there are 
 there are $rm-r$ 
choices for the distinguished
element $j_{r+1}$. Hence
\begin{align*}
        \left\lvert{\rm Cyc}_r(rm)\right\rvert & = (rm-1)_{r-1} \left\lvert{\rm Cyc}_r(rm-r)\right\rvert \\[6pt]
        & \qquad {}+ (rm-1)_{r-1}  (rm-r) \left\lvert{\rm Cyc}_r(rm-r)\right\rvert, 
    \end{align*} 
    which gives \eqref{r-2}. 
\qed

As per the recurrence
relations \eqref{r-1} and \eqref{r-2}, one can derive
the formulas for $\left\lvert{\rm Reg}_r(rm)\right\rvert$ and $\left\lvert{\rm Cyc}_r(rm)\right\rvert$, which result in the stronger version of the BMW inequality, 
i.e., \eqref{BMW-new}. Thus, for a prime $r\ge 3$, we obtain another combinatorial explanation of the monotone property. 

As noted in \cite{BMW-2000}, the case $r=2$ requires a stronger inequality, which they justify using generating functions.

\begin{lemma}
    For $m \ge 4$, we have
     \begin{align}
     \label{BMW-2}
    2\left\lvert{\rm Cyc}_4(4m)\right\rvert < \left\lvert{\rm Reg}_2(4m)\right\rvert.
    \end{align}
\end{lemma}

\proof
    For $m \ge 1$, by Lemmas \ref{recurrence-reg} and \ref{recurrence-cyc}, we obtain that
    \begin{align}
        &\left\lvert{\rm Reg}_2(4m)\right\rvert = (4m-1)^2\,(4m-3)^2 \left\lvert{\rm Reg}_2(4m-4)\right\rvert, \label{r-3} \\[6pt]
        &\left\lvert{\rm Cyc}_4(4m)\right\rvert = (4m-1)\,(4m-2)\,(4m-3)^2 \left\lvert{\rm Cyc}_4(4m-4)\right\rvert. \label{r-4}
    \end{align}
    In fact, the proofs of Lemmas \ref{recurrence-reg} and \ref{recurrence-cyc} 
    reveal that there is a bijection from ${\rm Reg}_2(4m)$ 
    to $[4m-1]^2\times[4m-3]^2\times{\rm Reg}_2(4m-4)$, 
    and there is also a bijection from ${\rm Cyc}_4(4m)$ 
    to $[4m-1] \times[4m-2]\times[4m-3]^2 \times {\rm Cyc}_4(4m-4)$.
    Clearly, the coefficient in \eqref{r-3}
    is greater than that in \eqref{r-4}, and it is just 
    a matter of formality to make this comparison in combinatorial
    terms. 
    Consequently, if   
    \begin{align*}
        2\left\lvert{\rm Cyc}_4(4m)\right\rvert < \left\lvert{\rm Reg}_2(4m)\right\rvert 
    \end{align*}
   holds for some value $m_0$, then it holds for
    all $m\geq m_0$. 
It is easily verified that we can choose $m_0=4$, since from \eqref{r-3} and \eqref{r-4}, we have
\[
 \frac{\left\lvert{\rm Reg}_2(4m_0)\right\rvert}{ \left\lvert{\rm Cyc}_4(4m_0) \right\rvert}=\frac{15}{14}\cdot\frac{11}{10}\cdot\frac{7}{6}\cdot\frac{3}{2}\cdot\frac{\left\lvert{\rm Reg}_4(0)\right\rvert}{ \left\lvert{\rm Cyc}_2(0)\right\rvert}=\frac{33}{16}>2.
\]
So the Lemma is proved.
\qed

In view of the above argument,  inequality \eqref{BMW-2} admits a combinatorial interpretation. 
Appealing to this inequality, we deduce that $p_2(n+1) < p_2(n)$ 
whenever $n+1 \equiv  0 \pmod {4}$, with the only exceptions
for $n=3,7,11$, see \cite{BMW-2000}. For these three special cases, 
we can look up the data in \cite{BMW-2000} or the sequence A247005 in the OEIS \cite{OEIS}. 
The values of 
$p_2(n)$ for $n=3,4,7,8,11,12$ are given as follows
\begin{align*}
    \frac{1}{2}, \,\frac{1}{2}, \, \frac{3}{8}, \,\frac{17}{48}, \,\frac{29}{96}, \,\frac{209}{720}.
\end{align*} 
Thus for $n=3,7,11$, 
the inequality $p_2(n+1) \le p_2(n)$ is valid with
equality attained only when $n=3$. Therefore,
for any prime $r$, a full combinatorial analysis is achieved. 

\section{The monotone property for prime powers}

Whereas the monotone property of $p_r(n)$ does not hold in general, numerical evidence raises hopes that it might hold for prime powers, as indicated in Table \ref{tb-2} for $r=4, 8, 9$
and $1\leq n \leq 12$.
\begin{table}[ht]
        \setlength{\tabcolsep}{8pt}
        \renewcommand{\arraystretch}{1.2}
	\centering
	\begin{tabular}{|c|cccccccccccc|}
		\hline
         \diagbox[width=2.8em, height=2.5em]{$r$}{\raisebox{4pt}{$n$}}
            & \raisebox{-1pt}{$1$} & \raisebox{-1pt}{$2$} & \raisebox{-1pt}{$3$} & \raisebox{-1pt}{$4$} & \raisebox{-1pt}{$5$} & \raisebox{-1pt}{$6$} & \raisebox{-1pt}{$7$} & \raisebox{-1pt}{$8$} & \raisebox{-1pt}{$9$} &\raisebox{-1pt}{$10$} & \raisebox{-1pt}{$11$} & \raisebox{-1pt}{$12$} \\	\hline
		$4$ & ${ 1}$ & ${ \dfrac{1}{2}}$ & ${ \dfrac{1}{2}}$ & ${\dfrac{3}{8}}$ & ${ \dfrac{3}{8}}$ & ${ \dfrac{5}{16}}$ & ${ \dfrac{5}{16}}$ & ${ \dfrac{53}{192}}$ & ${ \dfrac{53}{192}}$ & ${ \dfrac{95}{384}}$ & ${ \dfrac{95}{384}}$ & ${ \dfrac{29}{128}}$ \raisebox{23pt}{}\\[10pt]	\hline
		$8$ &  ${ 1}$ & ${ \dfrac{1}{2}}$ & ${ \dfrac{1}{2}}$ & ${ \dfrac{3}{8}}$ & ${ \dfrac{3}{8}}$ & ${ \dfrac{5}{16}}$ & ${ \dfrac{5}{16}}$ & ${ \dfrac{35}{128}}$ & ${ \dfrac{35}{128}}$ & ${ \dfrac{63}{256}}$ & ${ \dfrac{63}{256}}$ & ${ \dfrac{231}{1024}}$ \raisebox{23pt}{}\\[10pt]		\hline
        $9$ & ${ 1}$ & ${ 1}$ & ${ \dfrac{2}{3}}$ & ${ \dfrac{2}{3}}$ & ${ \dfrac{2}{3}}$ & ${ \dfrac{5}{9}}$ & ${ \dfrac{5}{9}}$ & ${ \dfrac{5}{9}}$ & ${ \dfrac{40}{81}}$ & ${ \dfrac{40}{81}}$ & ${ \dfrac{40}{81}}$ & ${ \dfrac{110}{243}}$  \raisebox{23pt}{} \\[10pt]		\hline
	\end{tabular}
        \caption{The values of $p_r(n)$.}
        \label{tb-2}
\end{table}

The main result of this section is as follows. 

\begin{theorem}
\label{prn-primepower}
    Let $n$ be a positive integer and $r=q^l$, 
    where $q$ is a prime, and $l \ge 1$, then  
    $p_r(n) \ge p_r(n+1)$.
\end{theorem}

Like the case for primes, this monotone property
stands on the following cases subject to
modulo conditions on $n+1$. First, we recall an
equality concerning $p_r(n)$ when $r$ is a prime power. Chernoff \cite{Chernoff-1994} 
established the following equality, and 
Lea\~nos, Moreno and Rivera-Mart\'inez \cite{Leanos-etal-2012} presented 
two proofs, with one using generating functions and the other being combinatorial. 

\begin{theorem}
\label{root-probability-power-1}
Let $q$ be a prime and $r=q^l$, $l \ge 1$. If $n+1 \not\equiv 0 \pmod q$, then   $p_r(n)=p_r(n+1)$.
\end{theorem}

For the remaining cases, we obtain the following relations. 

{\begin{theorem}
\label{root-probability-power}
Let $q$ be a prime, and $r=q^l$, $l \ge 1$.
\begin{enumerate}
    \item[(i) ] If $n+1 \equiv 0 \pmod {q}$ but 
    $n+1 \not\equiv 0 \pmod {qr}$, 
    then 
    \begin{equation} \label{PR1}
        p_r(n) \geq \frac{n+1}{n}p_r(n+1),
    \end{equation}
    with equality only when $n+1=kq$, where $k=1,2,\ldots,r-1$.
    \item[(ii) ] If $n+1 \equiv 0 \pmod {qr}$, 
    then $p_r(n) \ge p_r(n+1)$ with equality only when $r=2$ and $n=3$.
\end{enumerate}
\end{theorem}

To prove the above theorem, it is helpful to prepare 
some auxiliary inequalities. 
Even though these estimates can be
considerably improved, we will be content with the
coarse lower bounds in order to keep the proofs brief. 
Recall that for any $r$,
permutations with 
  an $r$-th root can be characterized in terms of the cycle type by Knopfmacher and Warlimont, see Wilf \cite[p. 158]{Wilf-2005}. In particular, we
 need the following criterion
when $r$ is a prime power. 

\begin{proposition}
\label{roots-cycles}
    If $r=q^l$ with $q$ being a prime number and $l \ge 1$, then a 
    permutation has an $r$-th root if and only if for any integer $i$, the number of cycles of length $iq$ is
 a multiple of $r$.
\end{proposition}

Next, we show that the two inequalities in
Lemma 3.2 and Lemma 3.3 in \cite{BMW-2000} for a
prime  $r$ can be extended
to a prime power. 
With the common notation $S_n^r$
for the set of permutations of $[n]$ with an $r$-th root, by Proposition \ref{roots-cycles}, we have
for a prime power $r=q^l$,
\begin{align} \label{RS}
    {\rm Reg}_q(n) \subseteq S^r_n .
\end{align}

Let ${\rm Cyc}_{q,\,r}(n)$ denote the set of permutations such that each cycle length is a multiple of $q$ and each cycle length occurs a multiple of $r$ times. 
The following relation is an extension of 
Lemma 3.2 in \cite{BMW-2000}. 

\begin{lemma}
\label{lemma-3.2}
For any $m \ge 1$, let $r=q^l$, where $q \ge 2$ (not necessarily a prime) and $l \ge 1$, we have 
\begin{align}
    \label{eq-3.2}
    \frac{\left\lvert {\rm Cyc}_{qr}(mqr)\right\rvert}{\left\lvert {\rm Cyc}_{q,\,r}(mqr)\right\rvert } \ge (mq)^{r-1}.
\end{align}
\end{lemma}

\proof
Let $\pi \in {\rm Cyc}_{q,\,r}(mqr)$. 
By definition, we assume that $\pi$ contains $k_ir$ cycles of length $iq $, 
where $k_i \geq 0$. 
For each $i$ with $k_i \neq 0$, partition the cycles of length
$iq$ into $k_i$ 
classes with each class containing
$r$ cycles. For each class $F$ of $r$ cycles of length $iq$, we proceed to generate a cycle of length $iqr$ out of the
elements in $F$. 
Running over all such classes $F$, we 
obtain permutations 
in ${\rm Cyc}_{qr}(mqr)$. 

First, let $A_1, A_2, \ldots, A_{r}$ be the cycles in $F$, where every cycle has length $iq$, that is,
arrange the cycles in $F$ in any specific linear
order.  
To form a cycle of length $iqr$, we  represent 
$A_1$ with the minimum element at the beginning.
Then break the cycles $A_2, A_2, \ldots, A_r$ into linear orders by  
starting with any element. There are
$iq$ ways to break a cycle of length $iq$ into a linear order.
Assume that $A_2', A_3', \ldots, A_r'$ are in 
linear orders by breaking the cycles $A_2, A_3, \ldots, A_r$,
respectively. 
Now we can form a cycle of length  $iqr$ by
adjoining $A_2', A_3', \ldots, A_r'$ successively
at the end of $A_1$. 
Evidently, the cycles formed in this way are
all distinct, and there are $(iq)^{r-1}$ of them 
that can be generated in this manner.  

Taking into account all classes $F$, 
we may produce certain permutations in 
${\rm Cyc}_{qr}(mqr)$. The number of 
permutations generated this way
equals $\prod_i(iq )^{(r-1)k_i}$.
Moreover, the range of $i$ in 
$\prod_i(iq)^{(r-1)k_i}$  can be restricted to those such that
$k_i\geq 1$. Given that 
$q\geq 2$, for any $k_i \ge 1$, we have
$(iq)^{k_i} \geq i qk_i$ and \[ \prod_i iq{k_i} \ge \sum_i iqk_i.\] 
Thus we see that 
\[\prod_i(iq)^{(r-1)k_i} \ge \left(\sum_i iqk_i\right)^{r-1}=(mq)^{r-1}, \]
where we have used the relation 
\[ \sum_{i} iqk_i  = m q, \]
because 
\[ \sum_{i} iq k_i r = m qr. \]
So the Lemma is proved.  
\qed

The following lemma is an extension of Lemma 3.3 in \cite{BMW-2000}. 

\begin{lemma}
\label{lemma-3.3}
 Let $r=q^l$, where $q \ge 2$ (not necessarily a prime) and $l \ge 1$. 
 Then, for any integer $m \ge 1$, we have 
\begin{align} \label{eq-3.3}
    \frac{\left\lvert{\rm Reg}_{q}(mqr)\right\rvert}{\left\lvert{\rm Cyc}_{q,\,r}(mqr)\right\rvert} > (mq)^{r-1}.
\end{align}
\end{lemma}

\proof
By definition, we have
\begin{align*}
    {\rm Cyc}_{qr}(mqr)  \subset  {\rm Cyc}_{q}(mqr),
\end{align*}
and hence $\left\lvert{\rm Cyc}_{qr}(mqr) \right \rvert <  
\left\lvert{\rm Cyc}_{q}(mqr) \right\rvert$.
In light of the stronger version of the BMW inequality \eqref{BMW-new}, we see that
\begin{align}
    \label{eq-BMW-r}
    \frac{\left\lvert{\rm Reg}_{q}(mqr)\right\rvert}{\left\lvert{\rm Cyc}_{qr}(mqr)\right\rvert} = \frac{\left\lvert{\rm Reg}_{q}(mqr)\right\rvert}{\left\lvert{\rm Cyc}_{q}(mqr)\right\rvert} \cdot \frac{\left\lvert{\rm Cyc}_{q}(mqr)\right\rvert}{\left\lvert{\rm Cyc}_{qr}(mqr)\right\rvert} > 1.
\end{align}
Comparing with \eqref{eq-3.2} shows that
\begin{align*}
    \frac{\left\lvert{\rm Reg}_{q}(mqr)\right\rvert}{\left\lvert{\rm Cyc}_{q,\,r}(mqr)\right\rvert} =\frac{\left\lvert {\rm Cyc}_{qr}(mqr)\right\rvert}{\left\lvert {\rm Cyc}_{q,\,r}(mqr)\right\rvert } \cdot\frac{\left\lvert{\rm Reg}_{q}(mqr)\right\rvert}{\left\lvert{\rm Cyc}_{qr}(mqr)\right\rvert} >(mq)^{r-1}, 
\end{align*}
as required.
\qed

The following lower bound of 
$\left\lvert S_{mqr}^r \right\rvert$ will also be needed
in the proof of Theorem \ref{root-probability-power}.

\begin{lemma}
Let $r=q^l$ be a prime power
greater than $2$. For any  $m \ge 1$, we have
\begin{align}
    \label{eq-part3-2}
      \left\lvert S_{mqr}^r \right\rvert > mqr \left\lvert{\rm Cyc}_{q,\,r}(mqr)\right\rvert .
\end{align}
\end{lemma}

\proof
    For the conditions  stated in the lemma, we obtain 
\[r-1=q^{l}-1 \ge l+1,\] 
thus, $(mq)^{r-1} \ge mq^{l+1}$. 
Thanks to \eqref{eq-3.3}, we find that
\begin{align*} 
    \frac{\left\lvert{\rm Reg}_{q}(mqr)\right\rvert}{\left\lvert{\rm Cyc}_{q,\,r}(mqr)\right\rvert} > (mq)^{r-1} \ge mq^{l+1}=mqr.
\end{align*}
Recalling \eqref{RS}, 
we get
\begin{align*}
    \left\lvert S_{mqr}^r \right\rvert \ge \left\lvert {\rm Reg}_{q} (mqr)\right\rvert > mqr \left\lvert{\rm Cyc}_{q,\,r}(mqr)\right\rvert,
\end{align*}
as claimed.
\qed

For now, we still need one
more inequality, that is, Corollary 2.16 in \cite{BMW-2000}. 
Given a permutation $\sigma$, 
 we may 
partition its set of cycles into two classes. 
Let $R_q(\sigma)$ and $S_{q}(\sigma)$ denote the 
permutation consisting 
of the $q$-regular cycles and the 
permutation consisting 
of the $q$-singular cycles of $\sigma$, respectively, in lieu of 
$\sigma_{(\sim q)}$ 
and $\sigma_{(q)}$ as used in 
\cite{BMW-2000}. 
We refer to $R_q(\sigma)$ and $S_{q}(\sigma)$ as the $q$-regular part and $q$-singular part of $\sigma$, respectively.  

The cycle type $\rho$ of a
permutation is defined to be the multiset of its
cycle lengths, often denoted by $1^{k_1} 2^{k_2} \cdots n^{k_n}$, meaning that there are $k_i$ cycles of length $i$ for $1\leq i \leq n$. We write $\lvert \rho \rvert$ 
for the sum of the cycle lengths in $\rho$.

Let $S_{\rho, \, q}(n)$, in place of  
${\rm DIV}_{\rho,\,q}(n)$ as used in \cite{BMW-2000}, denote the set of permutations of $[n]$ whose $q$-singular part has cycle type  $\rho$. In particular, $S_{\emptyset, \, q}(n)$ is
the set of $q$-regular permutations, 
i.e., ${\rm Reg}_q(n)$. For instance, the permutation $\sigma=(1\ 2)\,(3\ 4)\,(5\ 9\ 7\ 8)\,(6\ 10\ 11\ 13)\,(12)$ belongs to $S_{2^2 4^2,\,2}(13)$.

The inequality in \cite{BMW-2000} can 
be restated as follows.

\begin{proposition}
\label{prop-2.16} Let $n \geq 1$, $q\geq 2$ (not necessarily a prime), 
and let $\rho$ be a cycle type with 
$\lvert \rho \rvert \leq n$.
    If $n+1$ is a multiple of $q$, then 
    \begin{align*}
         \left\lvert S_{\rho, \, q}(n)\right\rvert \geq  \frac{1}{n} \left\lvert S_{\rho, \, q}(n+1)\right\rvert,
    \end{align*}
   where   equality is attained 
   if and only if $\rho=\emptyset$.  
\end{proposition}

In the event that $\rho=\emptyset$, 
the equality states that if $n+1$ is a multiple of $q$, then
\[ n\left\lvert{\rm Reg}_q(n)\right\rvert=\left\lvert{\rm Reg}_q(n+1)\right\rvert, \]
which follows from the bijection between ${\rm Reg}_q(n) \times [n]$ and ${\rm Reg}_q(n+1)$ 
constructed by B\'ona, McLennan and White, see Lemma 2.6 in \cite{BMW-2000}. 

We are now ready to prove Theorem \ref{root-probability-power}. 
 
\noindent{\it Proof of Theorem \ref{root-probability-power}.} Let us first introduce a notation. Given positive integers $q,r$, we say
 that a cycle type $\rho$ is
 $(q,r)$-divisible, denoted by $(q,r) \mid \rho$, if every cycle length in $\rho$ is divisible by $q$,  and
 for any integer $i$, the number of cycles of length $iq$ is
 a multiple of $r$. We also set $(q,r) \mid \emptyset$ by convention. 
 
 Since $r$ is a prime power $q^l$, Proposition \ref{roots-cycles} implies that a
permutation of $[n]$ belongs to $S_{n}^r$ if and only 
if the cycle type of its $q$-singular part is $(q,r)$-divisible. 
So we have
\begin{align*}
S_{n}^r=\bigcup_{\substack{\lvert\rho\rvert \le n
    \\[2pt] (q,\,r)\,\mid\,\rho}} S_{\rho, \,q}(n) .
\end{align*}
Hence
\begin{align}\label{S-1}
    \left\lvert S_{n}^r \right\rvert=\sum_{\substack{\lvert\rho\rvert \le n
    \\[2pt] (q,\,r)\,\mid\,\rho}} \left\lvert S_{\rho, \,q}(n)\right\rvert .
\end{align}
Again, by Proposition \ref{roots-cycles}, a 
permutation of $[n+1]$ is in $S_{n+1}^r$ 
if and only if its $q$-singular cycle type 
is $(q,r)$-divisible, namely,
\begin{align*}
    S_{n+1}^r  = \bigcup_{\substack{\lvert\rho\rvert 
    \le n+1 \\[2pt](q,\,r)\,\mid\,\rho }} S_{\rho, \,q}(n+1) .
\end{align*}
Considering the range of $\rho$, we get
\begin{align} \label{S-2}
   \left\lvert S_{n+1}^r\right\rvert = 
   \sum_{\substack{\lvert\rho\rvert \le n 
   \\[2pt] (q,\,r)\,\mid\,\rho}} \left \lvert S_{\rho, \,q} (n+1) \right\rvert 
   + \sum_{\substack{\lvert\rho\rvert = n+1 \\[2pt] (q,\,r)\,\mid\,\rho}}
   \left\lvert  S_{\rho, \,q}(n+1)  \right\rvert.
\end{align} 
Concerning the terms in \eqref{S-1} and in the first sum in \eqref{S-2}, 
given any cycle type $\rho$  
with $\lvert \rho\rvert \le n$ and $(q,r) \mid \rho$, 
Proposition \ref{prop-2.16} asserts that
if $n+1$ is a multiple of $q$, then 
\begin{align}
\label{eq-BMW-2.16}
     \left\lvert S_{\rho, \, q}(n)\right\rvert \geq  \frac{1}{n} \left\lvert S_{\rho, \, q}(n+1)\right\rvert,
\end{align}
where equality is attained 
if and only if $\rho=\emptyset$. Therefore, 
\begin{align*}
    \lvert S_{n}^r \rvert 
    &=\sum_{\substack{\left\lvert \rho \right\rvert \le n \\[2pt] (q,\,r)\,\mid\,\rho}} \left\lvert 
    S_{\rho, \,q}(n) \right\rvert
    \\[6pt]
    &\ge \frac{1}{n} \sum_{\substack{\lvert\rho\rvert \le n \\[2pt] (q,\,r)\,\mid\,\rho}}
    \left\lvert S_{\rho, \, q}(n+1)\right\rvert 
    \\[6pt]
    &=\frac{1}{n} \, \Bigg(\left\lvert S_{n+1}^r \right\rvert 
    -  \sum_{\substack{\left\lvert \rho \right\rvert = n+1 \\[2pt] (q,\,r)\,\mid\,\rho}} \left\lvert S_{\rho, \,q} (n+1) \right\rvert \Bigg).
\end{align*}
Consequently,
\begin{align}
    \label{eq-prime-power}
     n\left\lvert S_{n}^r \right\rvert \ge \left\lvert S_{n+1}^r \right\rvert -  \sum_{\substack{\lvert\rho\rvert = n+1 \\[2pt](q,\,r)\,\mid\,\rho}} \left\lvert S_{\rho, \,q} (n+1) \right\rvert.
\end{align}

We now proceed to prove (i).
Assume that $n+1$ is a 
multiple of $q$  but not a multiple of $qr$. We claim that 
for a cycle type $\rho$ with $\left\lvert \rho \right\rvert = n+1$ and $(q,r)\mid\rho$, 
\begin{align}
\label{eq-part1}
    S_{\rho, \,q} (n+1) =\emptyset.
\end{align}
Suppose to the contrary that there exists a permutation 
in $S_{\rho, \,q} (n+1)$. Under the condition that
$\rho$ is $(q,r)$-divisible, we have $\left\lvert\rho \right\rvert$ 
is a multiple of $qr$, but we also have $\left\lvert\rho \right\rvert = n+1$, 
which   contradicts   the condition 
that $n+1$ is not a multiple of $qr$.   
Utilizing the property
\eqref{eq-part1} and the relation
\eqref{eq-prime-power}, we get
\begin{align*}
     n\left\lvert S_{n}^r \right\rvert \ge \left\lvert S_{n+1}^r \right\rvert, 
\end{align*}
which is equivalent to \eqref{PR1}. This proves (i).  

To prove  (ii), assume that $n+1=mqr$.  
We shall proceed in the same fashion as the argument 
 given in \cite{BMW-2000} when $r$ is a prime. 
 The case $r=2$ has been taken care of in the 
 preceding section. 
So we may set our mind on the case when
$r$ is a prime power greater than $2$.  

In the notation $(q,r)\mid\rho$, we can write 
\begin{align} \label{CS}
    {\rm Cyc}_{q,\,r}(mqr) =\bigcup_{\substack{\lvert\rho\rvert = mqr \\ (q,\, r)\,\mid\,\rho}} S_{\rho, \,q} (mqr).
\end{align}
Inserting \eqref{CS} into \eqref{eq-prime-power} yields
\begin{align}
    \label{eq-part3-1}
     (mqr-1)\left\lvert S_{mqr-1}^r \right\rvert \ge \left\lvert S_{mqr}^r \right\rvert - \left\lvert{\rm Cyc}_{q,\,r}(mqr)\right\rvert.
\end{align}
Putting  \eqref{eq-part3-2} into \eqref{eq-part3-1} yields
\begin{align*}
     (mqr-1)\left\lvert S^r_{mqr-1} \right\rvert > \left(1-\frac{1}{mqr}\right)\left\lvert S^r_{mqr} \right\rvert.
\end{align*}
It follows that 
\begin{align*}
     mqr\left\lvert S^r_{mqr-1} \right\rvert > \left\lvert S^r_{mqr} \right\rvert.
\end{align*}
Thus we conclude that
$p_r(n) > p_r(n+1)$ when $n+1 \equiv 0 \pmod {qr}$. 
This proves (ii).

The conditions under which equality holds
 in (i) and (ii) are readily discerned, and so the theorem is proved. 
\qed

To conclude, we note that 
the monotonicity of $p_r(n)$ for prime powers $r$, as stated in Theorem \ref{prn-primepower}, is immediate
from Theorems \ref{root-probability-power-1} and \ref{root-probability-power}.

\vskip 6mm 
\noindent{\large\bf Acknowledgments.} 
We are deeply indebted to the referees for their close
scrutiny of the manuscript and for their constructive suggestions.

\end{document}